\documentclass[aps,prl,twocolumn,superscriptaddress,nobibnotes,floatfix,nofootinbib,longbibliography]{revtex4-2}

\usepackage[T1]{fontenc}
\usepackage{lmodern}
\usepackage{amsmath,amssymb,amsthm,bm}
\usepackage{mathtools}
\usepackage{xcolor}
\usepackage{microtype}
\usepackage[hidelinks]{hyperref}
\usepackage{url}
\usepackage{graphicx}
\usepackage[caption=false]{subfig}

\hypersetup{
  pdftitle={Quantifying the remaining search space in analogue optimisation},
  pdfauthor={Jona Nagerl, Younes Naceur, and Natalia G. Berloff}
}

\newcommand{\RR}{\mathbb{R}}
\newcommand{\TT}{\mathbb{T}}

\newcommand{\E}{\mathbb{E}}

\newcommand{\cE}{\mathcal{E}}

\newcommand{\gs}{\mathrm{gs}}

\begin{document}

\title{Quantifying the remaining search space in analogue optimisation}
\author{Jona Nagerl}
\email{jjn42@cam.ac.uk}
\affiliation{Department of Applied Mathematics and Theoretical Physics,
University of Cambridge, Wilberforce Road, Cambridge CB3 0WA, United Kingdom}
\author{Younes Naceur}
\email{Younes.Naceur@icfo.eu}
\affiliation{ICFO-Institut de Ciencies Fotoniques, 08860 Castelldefels (Barcelona), Spain}
\affiliation{LAAS-CNRS - Laboratoire d’Analyse et d’Architecture des Systèmes, Centre National de la Recherche Scientifique, Toulouse, France}
\author{Natalia G. Berloff}
\affiliation{Department of Applied Mathematics and Theoretical Physics,
University of Cambridge, Wilberforce Road, Cambridge CB3 0WA, United Kingdom}
\date{\today}

\begin{abstract}
Benchmarking analogue optimisers requires certifying not only how close a
returned energy is to the optimum, but also how many better configurations may
remain. Semidefinite relaxations provide rigorous bounds on the ground-state
energy and selected observables, yet reveal little about the size of this
residual search space. We introduce a deterministic, instance-specific upper
bound on the number of binary configurations whose energies improve upon a
returned state by a prescribed margin. Combined with a semidefinite energy
interval, this yields a complementary certificate of both the magnitude and
multiplicity of any possible improvement. To judge an analogue optimiser's native performance, we provide a generalisation of the 
framework to obtain bounds on the Haar volume of the remaining continuous phase-space. To further characterise the remaining search space, we extend our bounds to resolve regions in configuration space in both the binary and continuous setting. Our approach couples a
variational-Chernoff reduction with an elimination bound of the resulting partition function. Comparisons with exact enumeration and independent Haar sampling are consistent with the bounds on small instances. In numerical experiments at $N=1000$, even the largest family-median
bound restricts the potentially improving binary configurations to
less than $10^{-144}$ of the full configuration space, while the
largest upper bound on the corresponding Haar volume is of order
$10^{-499}$.
\end{abstract}

\maketitle

%=====================================================================
\section{Introduction}
\label{sec:introduction}

Many binary optimisation problems, including Max-Cut and Max-2-SAT, can be
mapped to the continuous dynamics of oscillator networks, coherent Ising machines, annealers, and other physical optimisers \cite{Yamamoto2017,Inagaki2016,McMahon2016,
WangRoychowdhury2019,KalininBerloff2018,Goto2021,Mohseni2022}. Finding the solution to the problem translates to the global minimisers of 
\begin{equation}
 H(x)=-\sum_{i<j}J_{ij}x_i x_j-\sum_i h_i x_i,
 \qquad x_i\in\{\pm1\}.
\label{eq:ising}
\end{equation}
Here $N$ is the number of variables, the coupling matrix $J\in\RR^{N\times N}$ is real symmetric
with zero diagonal, and $h\in\RR^N$ are the local fields. The global minimisers of $H$ encode optimal solutions to the original problem, while nonglobal local minima represent metastable candidate
solutions \cite{Lucas2014,BorosHammer2002}.

However, in the relevant regime, exactly computing the ground-state energy against
which the returned state should be judged, often is computationally intractable. Most benchmarking protocols compare attained objectives, success frequencies,
or time to a target \cite{King2015,Hamerly2019,BernalNeira2025}.
These quantities are indispensable for comparing devices, but none certifies
how far a particular output is from the unknown optimum.  

Baccari \emph{et al.} addressed this logical gap for quantum optimisers by computing certified lower bounds on the ground-state energy \cite{Baccari2020}.  Combined with the energy $U$ of a returned configuration, such a relaxation gives the rigorous interval
\begin{equation}
                 L\le E_{\gs}\le U .
\label{eq:energy-interval-intro}
\end{equation}
The interval quantifies how much improvement is energetically possible, but not whether potential improvements occupy half the state space or one configuration in $2^N$.

Benchmarking a solution by how many solutions it beats is an established idea in combinatorial optimisation. Domination analysis rates a heuristic by the number of solutions it is guaranteed to be better than \cite{GutinEtAl2003,GutinYeo2005,BerendSkienaTwitto2007}. For Max-Cut and Max-$k$-SAT, polynomial algorithms with a constant domination ratio exist \cite{AlonGutinKrivelevich2004}, and related results cover Boolean quadratic programs \cite{PunnenEtAl2013}. Instance-specific dominance certificates have also been constructed for heuristic outputs of TSP and MaxSAT instances \cite{BerendSkienaTwitto2018}. Motivated by this viewpoint, for a
returned binary state $x^s\in \{\pm1 \}^N$ with energy $E= H(x^s)$ we consider
\begin{equation}
 C_{\rm B}(E,\Delta)
 =
 \#\left\{
 x\in\{\pm1\}^N:
 H(x)\le E-\Delta
 \right\}.
\label{eq:headline-count}
\end{equation}
For a fixed returned state $x^s$, we use the shorthand
$C(\Delta):=C_{\rm B}(H(x^s),\Delta)$.Thus $C(\Delta)$ counts all configurations that match or improve upon $x^s$ by at least $\Delta$. For $\Delta=0$ the count includes $x^s$ and any degenerate configuration, in particular $-x^s$ if $h_i=0$. The case $\Delta > 0$ is relevant if one is interested in the number of solutions offering a substantial improvement. Since
\begin{equation}
 D(x^s):=\#\{x:H(x)\ge H(x^s)\}\ge 2^N-C(0),
\end{equation}
an upper bound on $C(0)$ gives a lower bound on the domination number $D(x^s)$ of the returned state $x^s$, that is the number of states with energy as least as large as the trial state.

Density-of-states methods estimate the population of each energy level \cite{WangLandau2001,ErmonGomesSelman2010} but do not give one-sided bounds. Hashing-based model counters give bounds that hold with a prescribed probability, but require many calls to an exact solver \cite{GomesSabharwalSelman2006,ErmonEtAl2013}. Jain \emph{et al.} derive deterministic upper bounds on the number of feasible solutions of binary integer programs \cite{JainEtAl2010}. An objective cutoff can be added to such a program as a constraint, but the method targets small structured programs rather than dense quadratic objectives. Related hardware work enumerates solutions under assumptions on sampling fairness and energy ordering \cite{Mizuno2024}. Another approach bounds satisfiability densities by concentration inequalities whose sensitivity constants come from an annealer, so the bound is valid only if the annealer finds the true maxima \cite{SahaiEtAl2020}. What we are interested in here: a deterministic upper bound on the number of configurations better than the returned state, that applies to the continuous phase states in which oscillator machines compute. We seek bounds that are valid for an arbitrary dense instance where for a fixed instance the count depends on $x^s$ only through $H(x^s)$, independently of how that state was found.

We derive a deterministic upper bound of this kind from variational Chernoff bounds \cite{RavikumarLafferty2004} combined with partition-function bounds \cite{JaakkolaJordan1996}. We express the improving fraction as an event probability under a uniform reference measure. Bounding this probability certifies the improving count independently of the optimiser's sampling or search procedure.

The remainder of the paper first introduces the phase-oscillator model
and semidefinite energy certificates. We then derive configuration-count
bounds using an elimination procedure related to that of
Jaakkola-Jordan \cite{JaakkolaJordan1996}, with curvature profiles reused
across energy thresholds, and extend the framework to continuous phases
while preserving pairwise elimination. We show how the resulting count
and volume bounds can be combined with energy intervals and observable
constraints, to further resolve the remaining search space. Finally, we assess the tightness of the bounds against exact
enumeration and independent Haar sampling, evaluate their scaling on
dense random instances and NP-hard encodings with up to $N=1000$
variables, and discuss limitations and extensions.

%=====================================================================
\section{Phase-based oscillator model and benchmark questions}
\label{sec:model}

\subsection{Phase dynamics and locking}

Following standard phase-based descriptions of oscillator Ising machines
\cite{WangRoychowdhury2019,ErementchoukShuklaMazumder2022,Steinerberger2023},
we model the analogue solver as a Kuramoto network with second-harmonic
injection locking. Related phase-reduced descriptions also arise in
gain-dissipative photonic networks, including coupled lasers and
nonequilibrium condensates
\cite{KalininBerloff2018,KalininBerloff2018Discrete}.

Oscillator $i$ is represented by a continuous phase $\theta_i\in\TT\colon=[0,2\pi)$.
We assume identical natural frequencies and work in a frame rotating at
their common frequency. For couplings $J_{ij}$, fields $h_i$, and
fixed locking strength $\kappa$, the deterministic dynamics are generated
by the Lyapunov function
\begin{equation} 
\begin{aligned}
\cE_\kappa(\theta)=
&-\sum_{i<j}J_{ij}\cos(\theta_i-\theta_j)\\
&-\sum_i h_i\cos\theta_i
-\frac{\kappa}{2}\sum_i\cos 2\theta_i .
\end{aligned}
\label{eq:phase-energy}
\end{equation}

The first term describes XY interactions, where positive $J_{ij}$ favours
equal phases, while negative $J_{ij}$ favours a phase difference of $\pi$.
The field $h_i$ couples oscillator $i$ to a fixed phase reference,
favouring $0$ for $h_i>0$ and $\pi$ for $h_i<0$.

The second-harmonic term connects the continuous XY model to binary
Ising spins. For $\kappa>0$, it favours $\theta_i\in\{0,\pi\}$,
which encode $x_i=+1$ and $x_i=-1$, respectively.
When $h=0$, this locking reduces the global $U(1)$ phase symmetry
to a twofold symmetry. At the exact binary configurations,
$\cos(\theta_i-\theta_j)=x_i x_j$, $\cos\theta_i=x_i$,
and $\cos 2\theta_i=1$, giving
\begin{equation}
\cE_\kappa(x)=H(x)-\frac{N\kappa}{2},
\label{eq:energy-connection-XY-Ising}
\end{equation}
for $H(x)$ from Eq.~\eqref{eq:ising}. Thus, locking reshapes the continuous energy landscape while leaving
the ordering of binary configurations unchanged. At finite locking
strength, we obtain a binary state by rounding each phase to the nearer
of $0$ and $\pi$, equivalently setting
$x_i=\operatorname{sign}(\cos\theta_i)$.

After rescaling the phase-relaxation time, we model the solver trajectories by deterministic gradient descent on the phase-energy landscape,
\begin{equation}
\begin{aligned}
\dot{\theta}_i
={}&-\frac{\partial\cE_{\kappa}(t)}{\partial\theta_i}\\
={}&-\sum_{j\ne i}J_{ij}\sin(\theta_i-\theta_j)
-h_i\sin\theta_i
-\kappa(t)\sin 2\theta_i,
\end{aligned}
\label{eq:phase-dynamics}
\end{equation}
where $\kappa(t)$ is the applied locking schedule. At fixed $\kappa$, the chain rule gives
\begin{equation}
\begin{aligned}
    \frac{\mathrm d\cE_\kappa}{\mathrm dt}
=\sum_i\frac{\partial\cE_\kappa}{\partial\theta_i}\dot{\theta}_i
=-\sum_i\left(\frac{\partial\cE_\kappa}{\partial\theta_i}\right)^2
\le 0,
\end{aligned}
\end{equation}
so the dynamics decrease the programmed phase energy. Here, experimental noise such as gain imbalance, detuning, and other
hardware effects would modify the trajectory, but not the certificates once
the realised $(J,h,\kappa)$ and returned state are specified.

Previous work reports that the best performance is obtained by
gradually increasing the locking strength $\kappa$
\cite{Farasat2026,FarasatShukla2026}. Following this strategy, our runs begin
with an unlocked XY exploration and then ramp $\kappa(t)$. If at the end of the run we look for a binary readout we apply $x_i=\operatorname{sign}(\cos\theta_i)$ to the final state $\theta$. This schedule allows the interaction energy to
relax before strong locking drives the phases towards binary states.

\subsection{Returned state and benchmark questions}

An analogue run returns a measured phase configuration
$\theta^s\in\TT^N$ with energy $E_\kappa=\cE_\kappa(\theta^s)$. We judge it by two questions: one about the
continuous landscape the machine actually navigates in, and one about the binary Ising problem encoded.

In the native setting, the dominance question we aim to answer is:
\emph{How much of phase space is better than the returned
state?} We measure this by the volume
\begin{equation}
 C_\kappa(E_\kappa, \Delta)=
 \mu_{\rm H}\!\left\{\theta\in\TT^N:
 \cE_\kappa(\theta)\le E_\kappa-\Delta\right\},
\label{eq:native-volume}
\end{equation}
where $\mu_{\rm H}$ is the normalised Haar measure on $\TT^N$ and
$\kappa$ is the locking strength at readout. This judges the machine in
its own landscape. Because $\cE_\kappa$ changes with $\kappa$, volumes
taken at different $\kappa$ are not comparable, so we report
$C_\kappa$ only at the realised $\kappa$. At $\kappa=0$ the landscape
is the XY objective itself.

On the other hand, we are interested in the number of logical states that improve upon the returned decoded solution. In that case we ask: \emph{How many different bit strings are better than the delivered answer?}
The phases are read out as
\begin{equation}
 x_i^s=\operatorname{sign}(\cos\theta_i^s),
\label{eq:readout}
\end{equation}
with $\operatorname{sign}(0)=+1$. The resulting bit string is assessed by
the count $C(\Delta)$ of Eq.~\eqref{eq:headline-count}, the number of
Ising configurations whose energy lies at least $\Delta$ below $H(x^s)$.
Once $x^s$, $J$ and $h$ are fixed, $C(\Delta)$ does not depend on
$\kappa$. The locking schedule affects it only through which bit string
the machine delivers.

The following sections derive rigorous upper bounds
$\widehat C(\Delta)$ on the binary count and
$\widehat C_\kappa(E_\kappa,\Delta)$ on the normalised continuous volume.

\section{Lower bound for the energy}
\label{sec:energy}

Before bounding the remaining search space, we obtain a
lower bound on the optimal energy which can be used complementarily for benchmarking solvers. Writing
$E_\kappa^\star=\min_{\theta\in\TT^N}\cE_\kappa(\theta)$,
the energy of a realised configuration provides an upper
bound $U_\kappa$ on $E_\kappa^\star$. A lower bound follows by expressing the quadratic
phase energy with a matrix variable. Let
$c_i=\cos\theta_i$, $s_i=\sin\theta_i$, and
$v=(1,c_1,\ldots,c_N,s_1,\ldots,s_N)^\top$.
The matrix $Y=vv^\top$ is positive semidefinite and satisfies
\begin{equation}
 Y_{00}=1,\qquad
 Y_{c_ic_i}+Y_{s_is_i}=1,\quad i=1,\ldots,N.
\end{equation}
Furthermore, the energy is linear in the entries of $Y$:
\begin{equation}
\begin{aligned}
 \cE_\kappa(Y)
 ={}&-\sum_{i<j}J_{ij}
       \bigl(Y_{c_ic_j}+Y_{s_is_j}\bigr)\\
 &-\sum_i h_iY_{0c_i}
 -\frac{\kappa}{2}\sum_i
       \bigl(Y_{c_ic_i}-Y_{s_is_i}\bigr).
\end{aligned}
\end{equation}
The optimisation problem for the optimal energy then reads
\begin{equation}
\begin{aligned}
 E_\kappa^\star=\min_Y\quad &\cE_\kappa(Y)\\
 \text{subject to}\quad
 &Y\succeq0,\qquad \operatorname{rank}Y=1,\\
 &Y_{00}=1,\\
 &Y_{c_ic_i}+Y_{s_is_i}=1,
 \quad i=1,\ldots,N.
\end{aligned}
\label{eq:direct-phase-sdp}
\end{equation}
Dropping the nonconvex rank-one constraint yields a
semidefinite programme (SDP) \cite{Laurent2009}, whose optimal value
$\ell_{\kappa,1}$ satisfies
$\ell_{\kappa,1}\le E_\kappa^\star$.
As the resulting matrices need not correspond to individual
phase configurations due to potentially larger ranks, this yields a lower bound to the original minimisation problem.

For binary optimisation problems, we instead write $v_{\rm B}=(1,x)^\top$
and minimise
\begin{equation}
 H(Y)=-\sum_{i<j}J_{ij}Y_{ij}-\sum_i h_iY_{0i}
\end{equation}
over $Y\succeq0$ with $Y_{ii}=1$ for $i=0,\ldots,N$.
The binary and phase PSD matrices therefore have dimensions
$N+1$ and $2N+1$, respectively, with $N+1$ scalar
normalisation constraints in either case.

These elementary relaxations coincide with the first
level of the Lasserre moment-sum-of-squares hierarchy
\cite{Lasserre2001,Laurent2009}. Higher relaxation orders introduce
higher-degree moments and additional
constraints, giving a systematic sequence of tighter
bounds that converges to the original solution in the binary and toric scenarios considered here. We use only the first level and construct its
matrix constraints explicitly in JuMP \cite{Lubin2023}, solving them with
MOSEK~\cite{MOSEK2025}.

The numerical solver output is converted into a certified
lower bound $L_{\kappa,1}\le\ell_{\kappa,1}$ by the certification procedure described in
Appendix~\ref{app:numerical-certification}.
Together with an outward-rounded upper evaluation of the
returned energy, this gives the certified interval
\begin{equation}
 L_{\kappa,1}\le E_\kappa^\star\le U_\kappa.
 \label{eq:phase-optimality-interval}
\end{equation}
The binary energy interval is obtained in the same way.

\section{Counting states below an energy threshold}
\label{sec:counting}

\subsection{Chernoff reduction}

We are now going to describe how variational Chernoff bounds on event probabilities \cite{RavikumarLafferty2004} can be used to obtain bounds on the number of states below a given energy threshold. We are interested in the number of states with energy at most $E$, which we denote by $C_{\mathrm{B}}(E) \equiv C_{\rm B}(E,0)$, which is defined in Eq.~\eqref{eq:headline-count}.

The idea is to introduce a probabilistic interpretation and describe $C_{\mathrm{B}}(E)$ as the expectation value of a random variable, and then to apply Markov's inequality to obtain a tractable quantity. This converts the counting problem into a log-partition-function
bound, for which tractable variational upper bounds have been
developed~\cite{WainwrightJaakkolaWillsky2005}.

Let $X$ be a random variable uniformly distributed  on $\{\pm1\}^N$.
For any parameter choice $\beta>0$, the probability of $H(X)\le E$ is equivalent to
$e^{-\beta H(X)}\ge e^{-\beta E}$. Applying Markov's inequality
to the nonnegative random variable $e^{-\beta H(X)}$ therefore
gives the exponential bound:
\begin{align}
 C_{\rm B}(E)
 &=2^N\mathbb{P}\!\left[H(X)\le E\right]\nonumber\\
 &\le 2^N e^{\beta E}
       \mathbb{E}_X\!\left[e^{-\beta H(X)}\right].
 \label{eq:chernoff-count}
\end{align}
The same bound holds trivially for $\beta=0$. Thus any upper bound
\begin{equation}
 B(\beta)\ge
 \log\mathbb{E}_X\!\left[e^{-\beta H(X)}\right]
\end{equation}
implies
\begin{equation}
 C_{\rm B}(E)\le
 \exp\!\left(N\log 2+\beta E+B(\beta)\right),
\end{equation}
Here $\log$ denotes the natural logarithm.

%This reduction from tail counting to a log-partition bound is classical \cite{Boucheron2013,DelaPenaGine1999} and underlies variational approaches to probabilistic inference \cite{RavikumarLafferty2004,WainwrightJaakkolaWillsky2005}.
While the above bound is valid for all $\beta$, its choice affects the tightness of the bound. Hence, in our implementation, we optimise over a given parameter range. If the exponent \(N\log 2+\beta E+B(\beta)\) is negative for any admissible $\beta$,
the corresponding upper bound on $C_{\rm B}(E)$ is strictly
less than one and hence certifies that no states with energy below $E$ exist.

\subsection{Recursively bounding the partition function}
For uniform binary spins, the exponential moment is the
normalised Ising partition function:
\begin{equation}
\begin{aligned}
 \mathbb{E}_X\!\left[e^{-\beta H(X)}\right]
 &=2^{-N}Z_{\rm B}(\beta),
 \\
 Z_{\rm B}(\beta)&=
 \sum_{x\in\{\pm1\}^N}e^{-\beta H(x)}.
 \end{aligned}
\end{equation}
Thus Eq.~\eqref{eq:chernoff-count} becomes
$C_{\rm B}(E)\le e^{\beta E}Z_{\rm B}(\beta)$.
Exact evaluation is intractable in
general~\cite{GoldbergJerrum2007,Nagy2024}, motivating
tractable upper bounds on the log-partition
function~\cite{WainwrightJaakkolaWillsky2005,WainwrightJordan2006}.
Several methods provide rigorous upper bounds on Ising partition functions,
including tree-reweighted and log-determinant relaxations
\cite{WainwrightJaakkolaWillsky2005,WainwrightJordan2006}, convex and
entropy-based bounds \cite{ElGhaouiGueye2008,BeylerBach2025}, and approximate
elimination schemes \cite{DechterRish2003,LiuIhler2011}. 

Any rigorous upper bound on the partition function can be inserted into the Chernoff construction to obtain a valid bound on the remaining configurations; the choice of bound affects the tightness and computational cost of the resulting certificate. Here we adopt the recursive quadratic bound of Jaakkola and Jordan \cite{JaakkolaJordan1996,JordanEtAl1999}. Its explicit recursive updates apply to dense coupling
matrices. The parametrisation developed below allows the
results of these updates to be reused when varying the
Chernoff parameter $\beta$ or the energy threshold $E$.
In particular, different returned configurations of the
same problem instance can be assessed without repeating
the recursive elimination for each configuration.
These properties motivate our implementation choice, without implying that this bound is uniformly tighter or faster than the alternatives.

The idea is to compute the expectation value by recursively eliminating single spins, summing over one spin while keeping the remaining spins fixed, and repeating this process.
Writing the contribution of spin $i$ to the energy as
$-x_i b_i$, with
\begin{equation}
 b_i=h_i+\sum_{j\ne i}J_{ij}x_j,
\end{equation}
the normalised sum over $x_i$ contributes
\begin{equation}
 \frac12\sum_{x_i=\pm1}e^{\beta x_i b_i}
 =\cosh(\beta b_i)
\end{equation}
to the partition function. Exactly eliminating one spin therefore adds $\log\cosh(\beta b_i)$
to the exponent, generally introducing interactions beyond
the quadratic Ising form. To preserve the original form of the Hamiltonian, to be able to work recursively, Jaakkola and Jordan bound this term by a quadratic envelope with parametrised curvature $c$,
\begin{equation}
\begin{aligned}
 \log\cosh z\le cz^2+d(c)\quad\text{with}\\
 \qquad
 d(c)=\sup_{z\in\RR}\{\log\cosh z-&cz^2\},
 \quad 0<c\le\frac12 .
 \end{aligned}
\label{eq:binary-envelope}
\end{equation}
We evaluate $d(c)$ by rigorously enclosing its maximiser
and the resulting supremum, with $d(1/2)=0$.

At each elimination step, the emerging $\cosh$ factor is bounded by a quadratic contribution. Defining $\gamma \coloneqq 2c \beta$ where $c$ is the curvature of the quadratic bound and
$\beta$ is the Chernoff parameter, the terms involving
the remaining spins are absorbed into the reduced Ising
model $\mathbb{S}^{N\times N} \rightarrow \mathbb{S}^{(N-1)\times (N-1)}$ through the updates :
\begin{align}
 J_{jk}&\leftarrow J_{jk}+\gamma v_jv_k,
 &&j<k,\quad j,k\ne i,\nonumber\\
 h_j&\leftarrow h_j+\gamma h_i v_j,
 &&j\ne i,
 \label{eq:rank-one-update}
\end{align}
where $v=J_{-i,i}$ contains the current couplings
to the spin being eliminated.

The bounded error shifts the bound at each elimination step by
\begin{equation}
 c\beta^2\bigl(h_i^2+\|v\|_2^2\bigr)+d(c).
\end{equation}

We use a common curvature $c$ across all elimination steps, which gives a fixed parameter $\gamma$, and eliminating $N-1$ spins produces the
contribution
\begin{equation}
 S(\gamma)=
 \sum_{\text{eliminated }i}
 \left(h_i^2+\|J_{-i,i}\|_2^2\right),
\end{equation}
where each term is evaluated immediately before eliminating
$i$. Together with the final spin's field $h_f(\gamma)$,
the quantity $S(\gamma)$ forms the stored \textit{profile}.
Summing over the final spin exactly yields the final partition function bound
\begin{align}
 B_{\rm B}(\gamma,c)
 ={}&\frac{\gamma^2}{4c}S(\gamma)
 +(N-1)d(c)\nonumber\\
 &+\log\cosh\!\left(
 \frac{\gamma h_f(\gamma)}{2c}
 \right).
 \label{eq:binary-profile-bound}
\end{align}
For every $\gamma>0$ and $0<c\le1/2$, this expression bounds
the logarithm of the normalised partition function at
$\beta=\gamma/(2c)$. The stored quantities $S(\gamma)$ and
$h_f(\gamma)$ depend only on $\gamma$; the remaining dependence
on $c$ is explicit.

Consequently,
\begin{equation}
 \log\mathbb{E}_X\!\left[e^{-\beta H(X)}\right]
 \le B_{\rm B}(\gamma,c),
 \qquad \beta=\frac{\gamma}{2c}.
 \label{eq:binary-mgf-bound}
\end{equation}
The elimination order affects tightness but not validity.
We greedily eliminate the spin with the smallest current
$h_i^2+\|J_{-i,i}\|_2^2$, minimising the immediate quadratic
contribution.

Numerically, for each $\gamma$ in a finite grid
$\Gamma\subset(0,\infty)$, we therefore perform one elimination
recursion and store the resulting profile. Once the profiles have been computed for all $\gamma\in\Gamma$,
each threshold $E$ requires only scalar optimisation over $c$
and comparison across the stored profiles.
For each profile, we minimise the full Chernoff exponent
$\gamma E/(2c)+B_{\rm B}(\gamma,c)$, rather than
$B_{\rm B}(\gamma,c)$ alone. Selecting the smallest result
across $\gamma\in\Gamma$ gives
\begin{equation}
\begin{aligned}
\small
 &C_{\rm B}(E)\le\\
 &\exp\!\left[
 N\log 2+
 \min_{\gamma\in\Gamma}
 \inf_{0<c\le1/2}
 \left\{
 \frac{\gamma E}{2c}+B_{\rm B}(\gamma,c)
 \right\}
 \right].
 \end{aligned}
\end{equation}
The optimal choice of $c$ and the selected $\gamma$ can
therefore depend on $E$, while the stored profiles are
reused unchanged.
The trivial bound $C_{\rm B}(E)\le2^N$ is also included.

Using a common curvature sacrifices the freedom to choose
a separate $c_i$ at every step, but makes the matrix recursion \eqref{eq:rank-one-update}
depend only on $\gamma$. With step-dependent curvatures,
changing an early $c_i$ generally changes all subsequent
couplings and requires repeating the recursion.
For $G=|\Gamma|$ values of $\gamma$, the implementation for a dense interaction graph
requires $O(GN^3)$ time and $O(N^2)$ working memory,
plus storage for the profiles. 

As the above discussion holds in exact arithmetic, numerical errors have to be accounted for when manipulating floating-point numbers. Appendix \ref{app:numerical-certification} describes this numerical error control while Appendix \ref{app:numerical-settings} specifies
the parameter grids, numerical searches, and measured runtimes
and memory usage for the numerical section.

\subsection{Extension to continuous phase space}
\label{sec:continuous-count}
We now generalise the discrete procedure to bound the volume of
continuous phase configurations with energy below $E_\kappa$.
We denote this normalised phase-space volume by
$C_\kappa(E_\kappa)\equiv C_\kappa(E_\kappa,0)$, as defined in
Eq.~\eqref{eq:native-volume}. Here, the uniform binary measure is
replaced by the normalised Haar measure $\mu_{\rm H}$ on the torus,
with $d\mu_{\rm H}(\theta)=(2\pi)^{-N}d\theta$. Thus,
$0\le C_\kappa(E_\kappa)\le1$, and the volume has the same
probabilistic interpretation as the normalised binary count.
To bound this volume, we extend the binary elimination procedure
to the corresponding continuous partition integral.

For $\beta>0$, the event $\cE_\kappa(\theta)\le E_\kappa$
is equivalent to
$e^{-\beta\cE_\kappa(\theta)}\ge e^{-\beta E_\kappa}$.
Applying Markov's inequality under $\mu_{\rm H}$ therefore
gives
\begin{equation}
\begin{aligned}
 C_\kappa(E_\kappa)
 % &\le e^{\beta E}M_{\rm T}(\beta),
 &\le e^{\beta E_\kappa}M_{\rm T}(\beta),
 \\
 M_{\rm T}(\beta)&=
 \int_{\TT^N}e^{-\beta\cE_\kappa(\theta)}
 d\mu_{\rm H}(\theta).
 \end{aligned}
\label{eq:torus-chernoff}
\end{equation}
The bound also holds for $\beta=0$.
Here $M_{\rm T}$ is the continuous partition integral with
respect to the normalised Haar measure, replacing the normalised
binary partition sum.
Our task is consequently to upper-bound $\log M_{\rm T}$.
As in the binary case, we integrate out one variable at a
time and bound the resulting factor by a quadratic envelope
that preserves the pairwise model.

For the current complex local field, write
\begin{equation}
 A_i=h_i+\sum_{j\ne i}J_{ij}e^{-\mathrm{i}\theta_j},
 \quad
 \beta A_i=re^{-\mathrm{i}\phi},
 \quad
 \lambda=\frac{\beta\kappa}{2}.
\end{equation}
Keeping the remaining phases fixed, exact integration over
phase $i$ produces the factor
\begin{equation}
 F_\lambda(r,\phi)=
 \int_0^{2\pi}
 e^{\lambda\cos(2\theta)+r\cos(\theta-\phi)}
 \frac{d\theta}{2\pi}.
 \label{eq:F-lambda}
\end{equation}
As in the binary case, the resulting marginal energy would change its structure by exactly eliminating a variable. To obtain a quadratic bound in $|A_i|$ that preserves the
pairwise structure during elimination, we proceed in two steps. We first bound the emerging
factor independently of the local-field orientation $\phi$.
Using the Fourier-Bessel expansion~\cite{DLMFBessel},
\begin{equation}
\begin{aligned}
\label{eq:fourier-bessel}
 F_\lambda(r,\phi)
 =&I_0(\lambda)I_0(r)\\ &+2\sum_{m\ge1}
 I_m(\lambda)I_{2m}(r)\cos(2m\phi),
 \end{aligned}
\end{equation}
where $I_m$ is the modified Bessel function of the first kind
of order $m$.
For $\lambda\ge0$, the coefficients are nonnegative, giving
the orientation-independent upper bound
\begin{equation}
 F_\lambda(r,\phi)\le F_\lambda(r,0)=:F_\lambda(r).
 \label{eq:F-orientation}
\end{equation}
Exact elimination adds $\log F_\lambda(r,\phi)$ to the
exponent, generally introducing higher-order interactions.
To preserve the pairwise structure, we hence further bound it by a
quadratic function of $r=\beta|A_i|$.
Define
\begin{equation}
 d_\lambda^{\rm T}(c)
 :=\sup_{r\ge0}\{\log F_\lambda(r)-cr^2\},
 \qquad c>0.
 \label{eq:torus-envelope}
\end{equation}
Then
\begin{equation}
 \log F_\lambda(r,\phi)
 \le cr^2+d_\lambda^{\rm T}(c).
 \label{eq:torus-envelope-bound}
\end{equation}
For each fixed pair $(\lambda,c)$, the envelope offset is
defined by a one-dimensional maximisation over $r\ge0$.
It is finite, as $F_\lambda(r)\le e^{\lambda+r}$ implies
$d_\lambda^{\rm T}(c)\le\lambda+1/(4c)$.
This maximisation is distinct from the subsequent minimisation
of the Chernoff exponent over $c$.

Writing $v=J_{-i,i}$, the squared local field satisfies
\begin{equation}
\begin{aligned}
\label{eq:complex-square}
 |A_i|^2
 ={}&h_i^2+\|v\|_2^2
 +2h_i\sum_jv_j\cos\theta_j\\
 &+2\sum_{j<k}v_jv_k\cos(\theta_j-\theta_k).
 \end{aligned}
\end{equation}
Hence the quadratic envelope adds only constant, field,
and pairwise terms when a spin is eliminated and preserves the structure recursively.
With a common curvature $c$ and $\gamma=2c\beta$, elimination
therefore gives the same coupling and field updates as
Eq.~\eqref{eq:rank-one-update}, with $d(c)$ replaced by
$d_\lambda^{\rm T}(c)$.

These updates depend only on $\gamma$.
We therefore choose a finite grid
$\Gamma\subset(0,\infty)$ and perform one elimination
recursion for each $\gamma\in\Gamma$.
For each recursion, we store the sum
$S(\gamma)=\sum_{\text{eliminated }i}
(h_i^2+\|J_{-i,i}\|_2^2)$, evaluated using the current
fields and couplings, and the final field $h_f(\gamma)$.
Together, these quantities define a profile that can be reused.

For a prescribed locking strength $\kappa\ge0$ and energy
threshold $E_\kappa$, we then optimise the Chernoff exponent
over $c>0$ at each fixed $\gamma$. In our numerical implementation, we evaluate a finite curvature grid and retain the smallest certified upper bound over all evaluated \((\gamma,c)\) pairs.
The parameters $\beta=\gamma/(2c)$ and
$\lambda=\gamma\kappa/(4c)$ are determined by $c$ and
are not optimised independently.

Applying the envelope to the first $N-1$ phases and
integrating the final phase exactly gives
\begin{align}
 B_{\rm T}(\gamma,c,\kappa)
 ={}&\frac{\gamma^2}{4c}S(\gamma)
 +(N-1)d_\lambda^{\rm T}(c)\nonumber\\
 &+\log F_\lambda\!\left(
       \frac{\gamma|h_f(\gamma)|}{2c}
     \right),
 \quad
 \lambda=\frac{\gamma\kappa}{4c},
 \label{eq:torus-profile-bound}
\end{align}
such that, for $\kappa\ge0$ and $\gamma,c>0$,
\begin{equation}
 \log M_{\rm T}\!\left(\frac{\gamma}{2c}\right)
 \le B_{\rm T}(\gamma,c,\kappa).
 \label{eq:continous-mgf-bound}
\end{equation}
Combining this with Eq.~\eqref{eq:torus-chernoff},
minimising the full Chernoff exponent over $c>0$,
and selecting the smallest result across $\gamma\in\Gamma$
yields
\begin{equation}
\begin{aligned}
 C_\kappa(E_\kappa)\le
 \min\!\Biggl\{1,\,
 \exp\!\Biggl[
 &\min_{\gamma\in\Gamma}\inf_{c>0}
 \biggl\{
 \gamma E_\kappa/2c\\
 &+B_{\rm T}(\gamma,c,\kappa)
 \biggr\}
 \Biggr]\Biggr\}.
\end{aligned}
 \label{eq:torus-volume-bound}
\end{equation}
The continuous envelope is well defined for every $c>0$.
Changing $c$ also changes
$\lambda=\gamma\kappa/(4c)$, so the scalar envelope must be
evaluated at the corresponding $(\lambda,c)$.
When $\kappa=0$, $F_0(r)=I_0(r)$; for $\kappa>0$, the
second-harmonic term enters through $F_\lambda$ and
$d_\lambda^{\rm T}(c)$.

Appendix \ref{app:numerical-certification} describes the certified
evaluation of $F_\lambda$ and $d_\lambda^{\rm T}$.
Our numerical implementation searches a finite curvature grid,
as specified together with the tolerances in
Appendix \ref{app:numerical-settings}.

\subsection{Tightness of the relaxations}
\label{sec:two-slacks}

Combining the Chernoff bound~\eqref{eq:chernoff-count} with
the elimination bound~\eqref{eq:binary-mgf-bound} (and similarly for ~\eqref{eq:continous-mgf-bound}) introduces
two distinct sources of looseness. To
separate them, let
\begin{equation}
 B_{\rm ex}(\beta)
 =\log\left[\E e^{-\beta H(X)} \right]
\end{equation}
denote the exact uniform log-partition function and define
\begin{align}
 \ell_Z(E)
 &=N+\frac{1}{\log 2}
   \inf_{\beta\ge0}
   \{\beta E+B_{\rm ex}(\beta)\},
 \nonumber\\
 \widehat\ell(E)
 &=N+\frac{1}{\log 2}
   \inf_{\beta\ge0}
   \{\beta E+B(\beta)\}.
\end{align}
If $C_{\rm B}(E)>0$ and
$\ell(E)=\log_2 C_{\rm B}(E)$, then
\begin{equation}
\begin{aligned}
 &\qquad\qquad\widehat\ell(E)-\ell(E)
 =\\&
 \underbrace{\ell_Z(E)-\ell(E)}
 _{\text{Chernoff relaxation}}
 +
 \underbrace{\widehat\ell(E)-\ell_Z(E)}
 _{\text{log-partition approximation}}
 \end{aligned}
\label{eq:slack-decomposition}
\end{equation}
Both terms are nonnegative. The first term is present even if the partition function is known exactly, so improving the partition-function bound $B$ can reduce only the second.

At intermediate percentiles there is no
reason for the Chernoff bound to be exact. Indeed, if
$E\ge\E \left [ H(X) \right ]$, convexity makes $\beta=0$ optimal, so even the exact partition function returns only the trivial bound $C_{\rm B}(E)\le2^N$.

In general, the Chernoff bound is well suited to the extreme low-energy tail. Under standard large-deviation regularity assumptions, optimisation of the exact log-partition function recovers the leading exponential decay rate of the tail (e.g. the low-energy part of the spectrum), although not generally its finite-size prefactor \cite{Chernoff1952,Boucheron2013}. At the binary optimum the bound becomes exact: if $E=E_{\rm gs}$, then
\begin{equation}
\begin{aligned}
 \lim_{\beta\to\infty}
 \left[
 N+\frac{\beta E_{\rm gs}+B_{\rm ex}(\beta)}{\log 2}
 \right]
 =\\
 \log_2\#\{x:H(x)=E_{\rm gs}\}.
 \end{aligned}
 \label{eq:groundstate_deg}
\end{equation}
Thus the Chernoff relaxation can dominate the error at moderate energy ranks, whereas the quality of the partition-function approximation becomes increasingly important in the deepest tail.

This focus on the low-energy tail also matches the practical role of the certificate. For typical or high-energy thresholds, the corresponding state fraction is not rare and can be estimated directly by uniform sampling, with
standard binomial confidence intervals
\cite{ClopperPearson1934}. If the target fraction is $p$, the relative sampling error scales as $(Mp)^{-1/2}$ for $M$ independent samples; direct sampling therefore becomes costly precisely when the better-than-threshold states are rare. More elaborate sequential Monte Carlo methods can address some rare-event regimes \cite{DelMoralDoucetJasra2006}, but require an additional sampling construction. We therefore use deterministic certificates for the rare, best-energy part of the landscape, while ordinary sampling remains available for its well-populated bulk.

The continuous phase-volume certificate admits the same
separation between Chernoff slack and partition-integral
slack, with normalised Haar measure replacing the uniform
binary measure. Its target is a volume rather than a
count: isolated continuous minima have zero Haar volume,
so the binary ground-state degeneracy identity from Eq.~\eqref{eq:groundstate_deg} does not
carry over. The tightness of the continuous certificate,
including the additional slack from recursive elimination,
is assessed numerically here.

\section{Observable certificates}
\label{sec:observables}

Besides the energy gap and the number of states, analysing the structure of an energy subspace could yield meaningful operational information. We're aiming to resolve the remaining search space along observables of the form
\begin{equation}
 O(x)=o_0+a^\top x+\frac12x^\top Kx,
\end{equation}
where $o_0\in\RR$, $a\in\RR^N$, and
$K=K^\top\in\RR^{N\times N}$. Useful observables could be magnetisation, cut imbalance, or the overlap $(x^s)^\top x$ with a trial state $x^s$. In particular, overlap bounds can reveal how many spins must change to obtain a better solution. For binary spins, the Hamming distance $d_H(x,x^s)$ between two configurations counts the differing spins and satisfies $d_H(x,x^s)=\frac12[N-(x^s)^\top x]$. Thus, ruling out improving solutions with overlap greater than $q$ proves that every improving solution must differ from $x^s$ in at least $(N-q)/2$ spins. To quantify the ranges of said observables within a given low-energy subspace, one is interested in computing 
\begin{align} \label{eq:observable_prob}
O_{\max}(E)=\max_{x}\quad &O(x)\\
\text{s.t.}\quad &H(x)\le E,\nonumber\\
&x_i^2=1,\quad i=1,\ldots,N.\nonumber
\end{align}
Following Wang \emph{et al.}~\cite{WangEtAl2024Observables},
we use a semidefinite relaxation with an energy upper bound
to bound an observable within a given energy sublevel
set. Writing $Y=(1,x)(1,x)^\top$, the energy and observable become
linear in the matrix entries:
\begin{equation}
\begin{aligned}
 H(Y)&=-\sum_{i<j}J_{ij}Y_{ij}-\sum_i h_iY_{0i},\\
 O(Y)&=o_0Y_{00}+\sum_i a_iY_{0i}
       +\frac12\sum_{i,j}K_{ij}Y_{ij}.
\end{aligned}
\end{equation}
The matrix satisfies $Y\succeq0$,
$Y_{ii}=1$, and $\operatorname{rank}Y=1$.
 Dropping the rank-one constraint  gives the SDP
\begin{equation}
\begin{aligned}
 \omega_1(E)=\max_Y\quad &O(Y)\\
 \text{subject to}\quad
 &Y\succeq0,\qquad Y_{ii}=1,\quad i=0,\ldots,N,\\
 &H(Y)\le E Y_{00}.
\end{aligned}
\label{eq:observable-sdp}
\end{equation}
Every binary configuration satisfying $H(x)\le E$
produces a feasible matrix, so
$O_{\max}(E)\le\omega_1(E)$.
This SDP is the first level of the corresponding
constrained moment-sum-of-squares hierarchy
\cite{Lasserre2001,Laurent2009}.
For continuous observables quadratic in
$(\cos\theta,\sin\theta)$, the same matrix rewrite gives
\begin{equation}
\begin{aligned}
 \omega_{\kappa,1}^{\rm T}(E)=\max_Y\quad &O(Y)\\
 \text{subject to}\quad
 &Y\succeq0,\qquad Y_{00}=1,\\
 &Y_{c_ic_i}+Y_{s_is_i}=1,\quad i=1,\ldots,N,\\
 &\cE_\kappa(Y)\le E Y_{00}.
\end{aligned}
\end{equation}
Every admissible phase configuration yields a feasible
matrix, so this SDP provides an upper bound to the observable within
the continuous energy sublevel set.

In our implementation we construct the observable SDP with JuMP \cite{Lubin2023} and solve it
using MOSEK \cite{MOSEK2025}.
The reported upper bounds are obtained by certifying
the numerical dual solution of the problem, as described in
Appendix \ref{app:numerical-certification}.

\subsection{Counting with an observable constraint}
Solving the constrained moment relaxation for Problem \eqref{eq:observable_prob} certifies that no state in the energy sublevel set exceeds the observable bound. However, it does not say how many states have a specified observable value or how that population changes as the threshold is varied. These questions are the primary purpose of the observable-resolved count we introduce here, first in the binary and then in the continuous case.

For a threshold value $o$, define the joint count
\begin{equation}
 C_O(E,o)=\#\{x:H(x)\le E,\ O(x)\ge o\}.
\label{eq:observable-count-def}
\end{equation}
For every $\beta,\mu\ge0$, as with the energy count, we apply the Chernoff bound with the uniform base measure and the statistics
$(H,O)$,
\begin{equation}
 C_O(E,o)\le2^N e^{\beta E-\mu o}
 \E\exp[-\beta H(X)+\mu O(X)].
\label{eq:joint-chernoff}
\end{equation}
The elimination remains applicable provided that the observable preserves the pairwise Ising form. 

The joint Chernoff bound can be interpreted as an exponential Lagrangian relaxation of the energy and observable constraints \cite{Touchette2009}. The nonnegative parameters $\beta$ and $\mu$ act as multipliers for the inequalities $H(x)\le E$ and $O(x)\ge o$, respectively. For $\beta>0$, the resulting weight can be written as
\begin{equation}
e^{-\beta H(x)+\mu O(x)}
=e^{-\beta[H(x)- t O(x)]},
\qquad t=\mu/\beta.
\end{equation}
This has the form of a generalised Gibbs ensemble, where $t$ represents the thermodynamic conjugate to the observable $O$ \cite{Jaynes1957}. In particular, when $O$ represents the particle number, $t$ plays the role of a chemical potential in the grand-canonical ensemble. For other observables such as magnetisation, it can be interpreted as an auxiliary field favouring configurations with larger observable values. 

For $\beta>0$, the tilted exponent
$-\beta H(x)+\mu O(x)$ remains quadratic.
Since $x_i^2=1$, its Ising parameters are obtained through
\begin{equation}
 J_{ij}\mapsto J_{ij}+\frac{\mu}{\beta}K_{ij}
 \quad(i<j),\qquad
 h\mapsto h+\frac{\mu}{\beta}a,
 \label{eq:observable-shift}
\end{equation}
retaining zero diagonal for $J$ and adding the constant
$\mu(o_0+\tfrac12\operatorname{tr}K)$ to the exponent. Several observables can be combined with nonnegative multipliers, at the cost of a higher-dimensional scalar search.

For binary states, the joint certificate also yields an
upper bound on the observable throughout the energy sublevel
set. Let $\widehat B_O(\beta,\mu)$ bound
$\log\E\exp[-\beta H(X)+\mu O(X)]$ from above.
Equation~\eqref{eq:joint-chernoff} then gives
\begin{equation}
 C_O(E,o)\le
 \exp\!\left[
 N\log 2+\beta E+\widehat B_O(\beta,\mu)-\mu o
 \right].
 \label{eq:observable-affine-bound}
\end{equation}
For $\mu>0$, the right-hand side is below one whenever
\[
 o>\frac{N\log 2+\beta E+\widehat B_O(\beta,\mu)}{\mu}.
\]
Since $C_O(E,o)$ is an integer, no configuration can
satisfy both constraints at such a threshold. Hence,
for a nonempty energy sublevel set,
\begin{equation}
 O_{\max}(E)\le
 \frac{N\log 2+\beta E+\widehat B_O(\beta,\mu)}{\mu}.
 \label{eq:observable-analytic-cutoff}
\end{equation}
We retain the smallest bound over the evaluated parameter
choices. Thus the observable bound follows directly from
the same count certificate, independently of any grid of
observable thresholds. Exclusion holds strictly above this
bound; equality alone does not certify emptiness.

For the overlap $O(x)=(x^s)^\top x$, the same count
certificate gives a minimum flip requirement.
Set $E_\Delta=H(x^s)-\Delta$, with $\Delta>0$, and define
\begin{equation}
\begin{aligned}
 r_{\rm cert}&(\Delta)=\\
 &\max\Bigl(\{0\}\cup
 \bigl\{r\in\{0,\ldots,N\}:\\
 &\qquad\qquad \widehat C_O(E_\Delta,N-2r)<1\bigr\}\Bigr).
 \end{aligned}
\end{equation}
Every state achieving the prescribed improvement must
then differ from $x^s$ in more than $r_{\rm cert}$ spins.
The baseline radius zero follows from $\Delta>0$;
$r_{\rm cert}=N$ excludes the improvement globally.

In the continuous scenario, for observables already present in the oscillator energy Eq.~\eqref{eq:phase-energy}, the unconstrained bound construction can be adapted similarly to the binary case. A phase magnetisation term $\sum_i a_i\cos\theta_i$ shifts $h$, a pair correlation shifts $J$, and the binarity observable
\begin{equation}
                  O_{\mathrm{bin}}(\theta)=\sum_i\cos2\theta_i
\end{equation}
simply replaces $\lambda$ by $\lambda+\mu$ in
Eqs.~\eqref{eq:F-lambda}-\eqref{eq:torus-profile-bound}. This gives a volume upper bound on phases that are simultaneously low in energy and highly locked. 

\section{Numerical evaluation}
\label{sec:numerics}

We evaluate the certificates against exact binary counts and independently
sampled continuous sublevel volumes, then apply the certification workflow
to simulated analogue-solver outputs at $N=250$-$1000$. Finally, we test
observable-resolved counts and compare the resulting observable bounds
with a dedicated SDP. 
The first-order moment relaxations for energy and observable
certification are constructed in JuMP and solved
with MOSEK~\cite{MOSEK2025}.

\paragraph{Certificate parameters.}
We evaluate a small, fixed grid of elimination parameters
adapted to the coefficient scale of each instance, and reuse
the resulting profiles across energy thresholds.
For each profile, we optimise over the curvature $c$, with
$\beta=\gamma/(2c)$ and additional scalar optimisation in the binary count case. Observable certificates additionally
use a fixed grid of tilts $t=\mu/\beta$.
We retain the smallest validated upper bound among the
evaluated candidates and the trivial bound.
Numerical error control is described in
Appendix \ref{app:numerical-certification}. The parameter grids and
search settings are specified in Appendix \ref{app:numerical-settings}.
\subsection{Reference checks and looseness}

We consider five instance families. The spin-glass controls are Gaussian Sherrington–Kirkpatrick (SK)
couplings $J_{ij}\sim\mathcal N(0,(N-1)^{-1})$, bimodal SK couplings
$J_{ij}=\pm(N-1)^{-1/2}$ with equal probability, and Gaussian SK couplings
with independent fields $h_i\sim\mathcal N(0,0.3^2)$
\cite{CrisantiRizzo2002}. The other two families are direct quadratic
encodings of weighted Max-Cut and weighted Max-2-SAT
\cite{GoemansWilliamson1995,Lucas2014}. The Max-Cut instances are complete
graphs with $w_{ij}=|g_{ij}|/\sqrt{N-1}$ for
$g_{ij}\sim\mathcal N(0,1)$. The Max-2-SAT instances contain $8N$ clauses
with uniformly selected distinct variables, independent literal signs, and
weights drawn uniformly from $[0.5,1.5]/\sqrt 8$.

For binary instances at $N=18$, exhaustive enumeration gives both
the exact sublevel fractions and the partition function. This
separates the Chernoff relaxation from the additional
quadratic-elimination error.

For continuous instances at $N=10$, a pilot sample of $10^6$
independent Haar configurations selects energy thresholds at prescribed
lower-tail fractions. A fresh sample of $2\times10^6$ configurations
estimates the sublevel volumes, with 95\% Clopper-Pearson confidence
intervals \cite{ClopperPearson1934}. To assess the error introduced by the Chernoff bound,
we also estimate
$ p_Z(E)=\inf_{\beta\ge0} e^{\beta E} \mathbb E_{\mu_{\rm H}} \!\left[e^{-\beta\mathcal E_\kappa(\theta)}\right]$.
We select $\beta$ using a separate pilot of $2^{19}$ configurations,
then evaluate the exponential moment at fixed $\beta$ using eight
independent batches of $2^{18}$ configurations. Variation across blocks estimates sampling uncertainty. Note that these Chernoff estimates are diagnostic and do not constitute rigorous bounds.

\begin{figure}[!t]
  \centering
  \subfloat[]{%
    \includegraphics[width=\columnwidth]
      {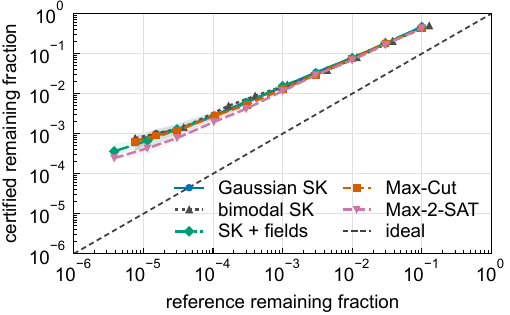}%
    \label{fig:small-calibration-binary}%
  }
  \par\vspace{-0.4em}
  \subfloat[]{%
    \includegraphics[width=\columnwidth]
      {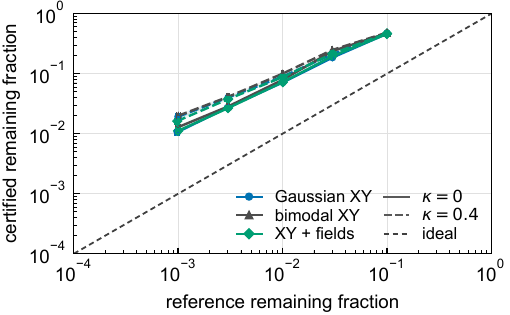}%
    \label{fig:small-calibration-phase}%
  }
  \caption{Calibration against independent references.
  (a) Exact binary fractions and certified fractions at $N=18$. For each
  family and exact-rank threshold, markers show the medians over ten seeds and
  bands show the interquartile range of the certified fraction.
  (b) Haar-volume estimates and certified phase fractions at $N=10$. Each curve displays the median over ten seeds. Markers show medians and bands show interquartile ranges.
  Horizontal error bars span the medians of the lower and upper endpoints of
  the per-instance 95\% binomial intervals obtained from $2\times10^6$ fresh
  samples. The diagonal denotes exact equality. Thresholds are fixed by exact
  rank in (a) and by an independent pilot in (b).}
  \label{fig:small-calibration}
\end{figure}

We observe no violations of the certified upper bounds in either
reference check [Fig.~\ref{fig:small-calibration}].
The certificates substantially restrict the state space, although
they remain loose relative to the reference fractions.
Table~\ref{tab:known-rank} separates the two sources of slack in
Eq.~\eqref{eq:slack-decomposition}. At the binary tenth-percentile
threshold, the Chernoff step dominates the looseness. At the ground-state threshold $E_{\gs}=\min_x H(x)$,
the Chernoff bound computed from the exact partition function is as expected  tight.
%as $\beta \to \infty$:
%\begin{equation}
%\begin{aligned}
% \inf_{\beta\ge0}\,2^{-N}e^{\beta E_{\gs}}Z_{\rm B}(\beta)
% &=2^{-N}\#\{x:H(x)=E_{\gs}\}.
%\end{aligned}
%\end{equation}

\begin{table}[t]
\caption{Median reference fractions and certified bounds.
Binary rows use five families and ten seeds per family
($N=18$, 50 instances); continuous rows use three native XY
ensembles and ten seeds per ensemble ($N=10$, 30 instances
per $\kappa$).
Binary reference fractions are exact, and binary Chernoff
values use the exact partition function.
Continuous reference fractions are Monte Carlo estimates,
marked $\dagger$; separate Chernoff estimates are not reported.
The certificate column contains deterministic upper bounds
obtained using quadratic elimination.}
\label{tab:known-rank}
\begin{ruledtabular}
% Dagger: Monte Carlo estimates; not rigorous bounds.
\begin{tabular}{lccc}
threshold & reference fraction & Chernoff & certificate\\
\hline
\multicolumn{4}{l}{Binary, $N=18$}\\
\hline
ground level & $7.63\times10^{-6}$ & $7.63\times10^{-6}$ & $4.35\times10^{-4}$\\
$0.1\%$ rank & $1.01\times10^{-3}$ & $6.80\times10^{-3}$ & $1.43\times10^{-2}$\\
$1\%$ rank & $1.00\times10^{-2}$ & $6.08\times10^{-2}$ & $7.72\times10^{-2}$\\
$10\%$ rank & $1.00\times10^{-1}$ & $4.32\times10^{-1}$ & $4.63\times10^{-1}$\\
\hline
\multicolumn{4}{l}{Continuous, $N=10$, $\kappa=0$}\\
\hline
$0.1\%$ target & $(9.91\times10^{-4})^{\dagger}$ & $(6.33\times10^{-3})^{\dagger}$ & $1.14\times10^{-2}$\\
$1\%$ target & $(1.00\times10^{-2})^{\dagger}$ & $(6.05\times10^{-2})^{\dagger}$ & $7.41\times10^{-2}$\\
$10\%$ target & $(1.00\times10^{-1})^{\dagger}$ & $(4.33\times10^{-1})^{\dagger}$ & $4.63\times10^{-1}$\\
\hline
\multicolumn{4}{l}{Continuous, $N=10$, $\kappa=0.4$}\\
\hline
$0.1\%$ target & $(9.96\times10^{-4})^{\dagger}$ & $(7.50\times10^{-3})^{\dagger}$ & $1.88\times10^{-2}$\\
$1\%$ target & $(1.00\times10^{-2})^{\dagger}$ & $(6.47\times10^{-2})^{\dagger}$ & $9.09\times10^{-2}$\\
$10\%$ target & $(1.00\times10^{-1})^{\dagger}$ & $(4.37\times10^{-1})^{\dagger}$ & $4.73\times10^{-1}$\\
\end{tabular}
\end{ruledtabular}
\end{table}

\subsection{Large instances}

We next assess the certificates on dense simulated instances with up to $N=1000$ variables, where we assess binary readouts and native continuous states separately.
In both cases, the analogue solver selects a candidate energy, and
the certificate bounds the state space below that energy. 
We consider system sizes $N\in\{250,500,750,1000\}$ and ten instance seeds per family. The binary benchmarks use the five
families defined above; continuous benchmarks comprise the three corresponding native XY ensembles: Gaussian couplings, bimodal couplings, and Gaussian couplings with
random fields. 

For each instance, the analogue solver evolves from 64 independent initial phase configurations
drawn uniformly from $\TT^N$. Each trajectory evolves deterministically
for dimensionless time 60 using 3000 forward-Euler steps of size
$0.02$. The locking protocol and candidate selection depend on the
benchmark, as described below.

\subsubsection{Binary readouts}

For binary optimisation, the locking strength remains zero for the
first 35\% of the evolution and then increases linearly to
$\kappa=1.2$. We round every final phase configuration and select the
readout with the lowest Ising energy.

At the returned energy $H(x^s)$, the binary certificate gives an
upper bound $\widehat C_{\rm B}$ on the number of configurations
satisfying $H(x)\le H(x^s)$. To compare different system sizes, we
report the binary count rate, $\widehat r_{\rm B}
 =\min\{1,\log_2\widehat C_{\rm B}/N\}$.
This normalisation accounts for the exponential growth of the
binary state space: a value $\widehat r_{\rm B}$ certifies at most
$2^{N\widehat r_{\rm B}}$ configurations, or a remaining fraction
of at most $2^{-N(1-\widehat r_{\rm B})}$. A rate of zero certifies that $x^s$ is the unique ground state, whereas a rate of one gives no improvement over the trivial bound of $2^N$ configurations.
Since \(x^s\) is itself counted, \(\widehat C_{\mathrm B}<2\) proves that \(x^s\) is the unique ground state; for \(h=0\), the count is even, and \(\widehat C_{\mathrm B}<4\) proves that \(x^s\) and \(-x^s\) are the only ground states.

\begin{figure}[!t]
  \centering
  \subfloat[]{%
    \includegraphics[width=.98\columnwidth]
      {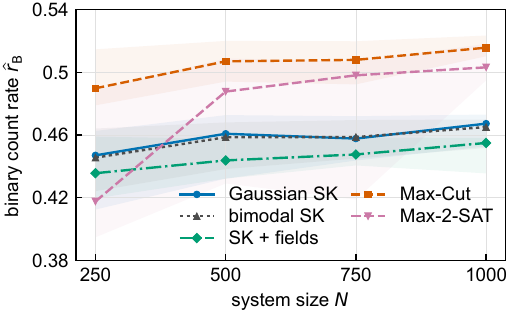}%
    \label{fig:dense-scaling-incumbent}%
  }
  \par\vspace{-0.4em}
  \subfloat[]{%
    \includegraphics[width=.98\columnwidth]
      {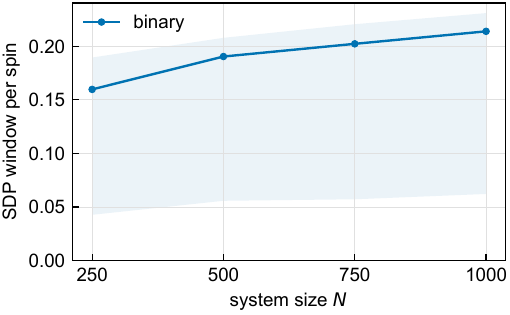}%
    \label{fig:dense-scaling-sdp-binary}%
  }
  \caption{Large binary benchmarks.
  (a) Binary count rate $\widehat r_{\rm B}$
  at the selected binary readout.
  Lines show medians over ten instance seeds per family and size;
  bands show the full seed range.
  (b) Gaps between the returned Ising energies and their SDP lower
  bounds, divided by $N$. The line shows the median over all five
  families and ten seeds per family, giving 50 instances at each
  size; the band shows the full range.}
  \label{fig:dense-scaling}
\end{figure}

At $N=1000$, the largest family-median value of $\widehat r_{\mathrm B}$ is about $0.52$ [Fig.~\ref{fig:dense-scaling}(a)]. For an instance with $\widehat r_{\mathrm B}\le0.52$, at most a fraction $2^{-0.48N}$ of all binary configurations match or improve on the returned state. Exponentially many such configurations may remain, and a small count rate does not by itself imply a small energy gap.

The SDP lower bound at relaxation order $d=1$ provides this interval. Its median width grows from $0.160$ per spin at $N=250$ to $0.214$ at $N=1000$ [Fig.~\ref{fig:dense-scaling}(b)]. The width $U-L=(U-E_{\mathrm{gs}})+(E_{\mathrm{gs}}-L)$ combines the suboptimality of the returned state with the slack of the relaxation. For zero-field Gaussian SK with our normalisation $J_{ij}\sim\mathcal N(0,(N-1)^{-1})$, the optimal first-order SDP value, denoted here by
$\ell_1$, satisfies $\ell_1/N\to-1$ \cite{MontanariSen2015}, while $E_{\mathrm{gs}}/N\to-P_*$ with $P_*\simeq0.7632$~\cite{DemboMontanariSen2015}. Hence, there is an analytical gap between the first order relaxation and the binary optimum which persists even if the analogue solver converged to the global optimum. Our variance $(N-1)^{-1}$ differs from the standard $N^{-1}$ convention by an energy scaling factor $\sqrt{N/(N-1)}\to1$, leaving both limits unchanged. Hence, a nonzero window per spin persists even for an optimal readout in the above scenario.

\paragraph{Improvement curves}

At $N=1000$, we vary the required improvement $\Delta$ for
all ten instances of each family, reusing the returned binary
states and certificate settings of the experiments above. 

\begin{figure}[!t]
\centering
\includegraphics[width=\columnwidth]{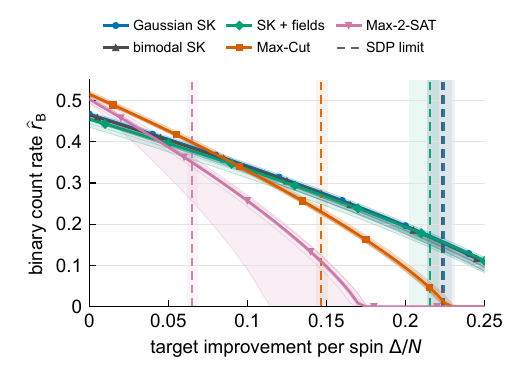}
\caption{Improvement curves at $N=1000$ for all five families,
plotted against the required improvement per spin $\Delta/N$.
Thin curves show individual instances; solid curves and shaded
bands show the median and full range over ten seeds.
Colour-matched dashed vertical lines mark the median of the SDP-derived improvement upper bounds
$\Delta_{\rm SDP}^{\rm UB}/N$, and vertical bands span its full range
over the ten seeds.
The curves show the count certificate separately from this
SDP exclusion, using the displayed rate
$\widehat r_{\rm B}(\Delta)
=\max\{0,\min\{1,\log_2\widehat C(\Delta)/N\}\}$.
Empty events are certified when the untruncated bound satisfies
$\widehat C(\Delta)<1$.}
\label{fig:improvement-curves}
\end{figure}

Figure~\ref{fig:improvement-curves} shows the binary count rate
$\widehat r_{\rm B}(\Delta)$ against the improvement per spin
$\Delta/N$, using 51 equally spaced thresholds spanning
$0\le\Delta/N\le0.25$. At zero improvement, the medians and
seed ranges reproduce the $N=1000$ large-instance results. We additionally provide the energy lower bound by the SDP. For each instance, the SDP energy lower bound $L_1$ gives an upper limit $\Delta_{\rm SDP}^{\rm UB}:=H(x^s)-L_1$ on the achievable
improvement: no configuration can satisfy $H(x)\le H(x^s)-\Delta$ when $\Delta>\Delta_{\rm SDP}^{\rm UB}$. The SDP lower bound identifies the improvement beyond which no configuration can exist, whereas the count certificate may still return a nonzero upper bound on the number of states in this region. Together, the two certificates bound both the number of configurations achieving a prescribed improvement and the maximum improvement that is possible.

\subsubsection{Native continuous states}

For continuous optimisation we keep $\kappa=0$ throughout the evolution and select the final configuration with the lowest XY energy $\mathcal E_0(\theta^s)$. The continuous certificate bounds the normalised Haar volume \eqref{eq:native-volume} at $\Delta=0$, that is, the probability that a uniformly sampled phase configuration matches or improves the returned energy. It is not an attraction-basin volume or a success probability of the optimiser. To compare system sizes, we report
$\widehat q_{\rm T}=\widehat p_{\rm T}^{\,1/N}$, with
\begin{equation}
 \widehat p_{\rm T}
 =\min\!\left\{1,
 \widehat C_0\!\left(\mathcal E_0(\theta^s),0\right)
 \right\}.
\end{equation}

\begin{figure}[!t]
  \centering
  \subfloat[]{%
    \includegraphics[width=.98\columnwidth]
      {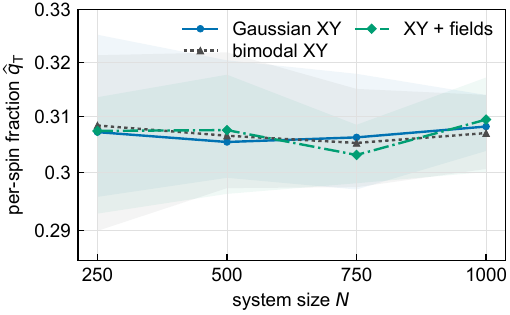}%
    \label{fig:xy-scaling-volume}%
  }
  \par\vspace{-0.4em}
  \subfloat[]{%
    \includegraphics[width=.98\columnwidth]
      {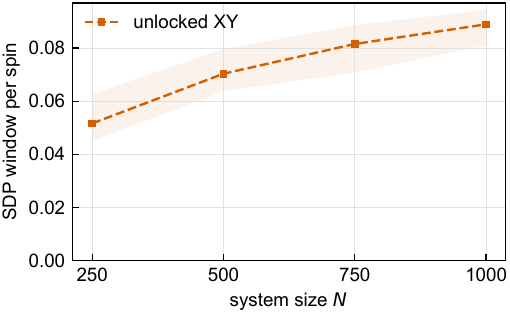}%
    \label{fig:xy-scaling-sdp}%
  }
  \caption{Large native continuous benchmarks, with $\kappa=0$.
  (a) Per-spin Haar-volume bound $\widehat q_{\rm T}$ at the
  returned XY energy. Lines show
  medians over ten seeds per ensemble and size; bands show the
  full seed range.
  (b) Gaps between the returned XY energies and their SDP lower
  bounds, divided by $N$. The line shows the median over three
  ensembles and ten seeds per ensemble, giving 30 instances at
  each size; the band shows the full range.}
  \label{fig:xy-scaling}
\end{figure}

Figure~\ref{fig:xy-scaling}(a) shows these bounds for the three native XY ensembles. Among the 30 instances at $N=1000$, the least restrictive bound is $\widehat q_{\rm T}\simeq0.317$, a Haar-volume bound of order $10^{-499}$. As in the binary case this does not show that the returned state is optimal. The SDP energy interval again complements the volume bound, with median width growing from $0.052$ per spin at $N=250$ to $0.089$ at $N=1000$ [Fig.~\ref{fig:xy-scaling}(b)].

\subsection{Observable-resolved counts}

Observable-resolved certificates bound how many states satisfy
both an energy cutoff and an observable constraint and can hence reveal additional information about the remaining search space. We test this extension using the signed magnetisation $M/N=N^{-1}\sum_i x_i$.

\begin{figure}[!htb]
  \centering
  \subfloat[]{%
    \includegraphics[width=\columnwidth]
      {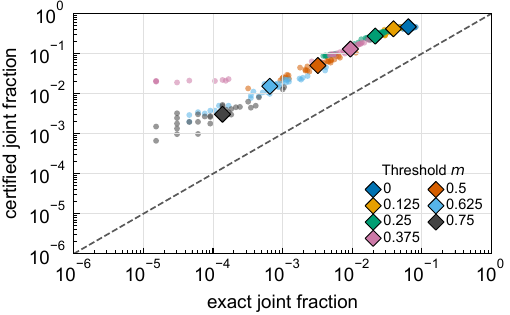}%
    \label{fig:observable-certificate-exact}%
  }
  \par\vspace{-0.4em}
  \subfloat[]{%
    \includegraphics[width=\columnwidth]
      {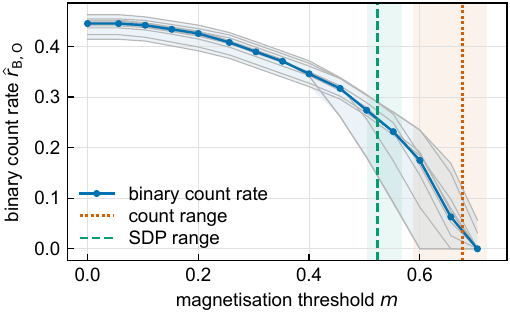}%
    \label{fig:observable-certificate-magnetisation}%
  }
\caption{Observable-resolved certificates.
  (a) Exact versus certified joint fractions at $N=16$, with each instance's
  energy cutoff fixed at its tenth energy percentile. Both axes show
  fractions of all binary configurations satisfying $H\leq E$ and
  $M/N\geq m$. Colours identify the signed magnetisation thresholds listed
  in the legend. Circles represent individual instances;
  diamonds give the coordinate-wise medians of the same-colour points.
  The dashed line denotes equality.
  (b) Certified magnetisation profile at $N=250$ for Gaussian SK. Thin curves
  show individual instances, and the blue line and band show the median and
  full range over ten instances. Vertical lines show the median count-derived
  and SDP range bounds, and the corresponding shaded vertical regions show
  their full ranges over the same ten instances.
  The count-derived range is obtained from
Eq.~\eqref{eq:observable-analytic-cutoff}.}
  \label{fig:observable-certificate}
\end{figure}

At $N=16$, we enumerate all binary configurations for ten instances
from each of the five families. For each instance, we fix the
energy cutoff $E$ at its tenth percentile and compare exact and
certified counts satisfying both $H\le E$ and the magnetisation bound $M/N\ge m$, where $m\in\{0,0.125,0.25,0.375,0.5,0.625,0.75\}$.
Dividing both counts by $2^{16}$ gives the exact and certified joint
fractions plotted in Fig.~\ref{fig:observable-certificate}(a). The resulting 350 instance-threshold comparisons include
41 cases with no configuration satisfying both constraints.
The certificate proves emptiness in 10 of these cases,
with $\widehat C_O(E,Nm)<1$, yielding an observable bound on the energy subspace.
All 41 zero-count cases are omitted from the logarithmic
axes, leaving 309 points in the figure.
Each point represents one instance-threshold pair, coloured by
$m$. For each colour, the diamond's horizontal coordinate is the median
exact fraction and its vertical coordinate is the median certified
bound across the displayed instances. These medians are taken separately
at fixed $m$. 
All points lie above the equality line, consistent with the bound on every tested instance.

Figure~\ref{fig:observable-certificate}(b) shows the
logarithmic joint count bound per spin,
\begin{equation}
 \widehat r_{\mathrm B,O}(m)
 = \max\!\left\{0,\min\!\left\{1,
 \log_2\widehat C_O(E,Nm)/N
 \right\}\right\},
 \label{eq:observable-count-rate}
\end{equation}
where $O(x)=M(x)=\sum_i x_i$ and the observable constraint
is $M(x)/N\ge m$.
Smaller values indicate tighter bounds on the number of
low-energy configurations with large magnetisation.
For display, negative rates are set to zero; empty events
are identified from the untruncated bound
$\widehat C_O(E,Nm)<1$.

At $N=250$, we evaluate the joint-count bounds at the
energy of the returned binary state. 
We obtain the count-derived magnetisation bound from
Eq.~\eqref{eq:observable-analytic-cutoff}, with
$O=M$, $o=Nm$, and $\mu=\beta t$.
We divide the resulting observable bound by $N$ and retain
the smallest upper endpoint over the evaluated parameter
choices with $t>0$, including the physical bound $M/N\le1$.
This uses the same candidates as the joint-count curves
but is independent of the plotted magnetisation grid.
The median count-derived bound over the ten instances is
$0.677$, compared with $0.524$ for the SDP.

Although the SDP gives a tighter magnetisation bound,
the joint-count curves additionally quantify how many
low-energy configurations may remain above each
magnetisation threshold. They therefore provide information
about the population of the energy sublevel set beyond
its extremal observable values.
%=====================================================================
\section{Discussion and outlook}
\label{sec:discussion}

The main contribution of this work is an end-to-end certificate that provides information about the remaining state-space in analogue optimisation. The SDP interval bounds how much
improvement remains energetically possible, while the count or volume curve
bounds how many states can realise each improvement. We develop tailored methods for bounding the remaining discrete logical space as well as the remaining native phase-space. To scale our procedure to experimentally relevant system sizes, we combine Chernoff reductions with controlled partition function bounds by Jaakkola and Jordan which we extend to the continuous scenario. To further characterise the remaining state space, we extend the above certificates to observable-resolved regions therein. 

Exact enumeration and
independent phase sampling validate the implementation, and the dense tests
show that both quantities remain easily computable up to $N=1000$. At this size the binary as well as the continuous certificate can certifiably exclude an exponentially large state space. The observable-resolved counts also give non-trivial statements about the remaining state space and, for binary states, can certify that no configurations satisfy the given constraints. The developed framework therefore offers a complementary work to extremal SDPs, which only certify an interval in which the optimum lives but do not offer access to the count of the states with certain qualities.

The practical relevance of our work is as a calibration layer around an optimiser. Each algorithmic or physical output can be mapped to
a certified energy interval and improvement curve. These paired quantities
can be used to compare integration times, noise levels, locking schedules, or
hardware configurations. They
complement success probability and time-to-target measurements by asking a
different question: after a state has been returned, how much of the search
space can still contain a prescribed improvement?

The numerical results also reveal limitations. At ordinary percentiles,
uniform sampling can estimate the state fraction directly and more tightly
with binomial confidence intervals \cite{ClopperPearson1934}; the deterministic
certificate is most useful in the rare low-energy tail, where such events are
difficult to observe. However, even there, an upper bound may still leave
exponentially many states. 

The next step is therefore to further tighten the certificate. Possible routes could be to improve the current bound especially improving the hyperparameter optimisation or to compare alternative rigorous
partition-function bounds.
The preferred choice may depend on the instance structure, the
energy range of interest, and the number of threshold queries.
Such comparisons could improve the trade-off between certificate
tightness and computational cost.
For problems mapped to an Ising Hamiltonian, constraint-aware elimination and conditioning hierarchies could tighten the reference space. Higher-order interactions and
multistate variables would extend the framework beyond quadratic binary
models \cite{Bybee2023,ChengLin2025Potts}.

\section*{Acknowledgements}
The authors thank Victor Magron for useful discussions. The authors acknowledge the use of GPT 6 Astra and Claude Opus 5.5 for assistance with brainstorming ideas, code development, and language improvement in the final manuscript. The final content, analysis, mathematical results, and conclusions remain the sole responsibility of the authors. 

J.N. acknowledges the support from the EPSRC through PhD studentship EP/W524633/1. 

Y.N. acknowledges support by European Union’s HORIZON–MSCA-2023-DN-JD programme under the Horizon Europe (HORIZON) Marie Skłodowska-Curie Actions, grant agreement 101120296 (TENORS), the Government of Spain (Severo Ochoa CEX2019-000910-S, FUNQIP and NextGeneration EU PRTR-C17.I1) and European
 Union (QSNP, 101114043 and PASQuanS2.1, 101113690), Fundació Cellex, Fundació Mir-Puig, and Generalitat de Catalunya (CERCA program).

N.G.B. acknowledges the support from  the HORIZON EIC Pathfinder Challenges project HEISINGBERG (grant 101114978), the Weizmann–UK Make Connection grant (grant 142568), and the EPSRC UK Multidisciplinary Centre for Neuromorphic Computing (grant UKRI982).

\section*{Code Availability}
The numerical implementation of the count and SDP certificates is publicly available under \url{https://github.com/jonaeg/StateSpaceCountCertificate}.

%=====================================================================
\appendix
\setcounter{secnumdepth}{1}

\section{Numerical certification}
\label{app:numerical-certification}

The discussion on how to obtain bounds on the state count, in the binary as well as in the continuous setting applies in exact arithmetic. When performing the algorithm on a computer using floating-point arithmetic, one hence has to take into account the finite-precision when obtaining rigorous bounds. 

Computed scalar floating-point quantities are enclosed rigorously using Arb~\cite{Johansson2017arb} through
\texttt{python-flint} with at least 160-bit precision.
Array calculations use double precision, with explicit
rounding-error bounds included in the certificates.
The SDP solver uses primal-feasibility, dual-feasibility,
and relative-gap tolerances of $10^{-8}$.
These tolerances control the numerical search, To certify numerical values we bound the residual errors
of each returned candidate, as described below.

\paragraph{Elimination errors.}
We first describe how to bound the errors which accumulate through the elimination bound to the partition function. At the step eliminating
spin $i$, let $J,h$ denote the current stored coefficients,
$v=J_{-i,i}$, and $\widehat J,\widehat h$ the stored
coefficients after the update. Define
\begin{equation}
\begin{aligned}
 \delta J^{(i)}_{jk}
 &=\widehat J_{jk}-(J_{jk}+\gamma v_jv_k),
 \\
 \delta h^{(i)}_j
 &=\widehat h_j-(h_j+\gamma h_i v_j),
\end{aligned}
\end{equation}
where the expressions in parentheses are evaluated
exactly from the current stored coefficients.
We compute an upper bound
\[
 \varepsilon_i\ge
 \sum_{\substack{j<k\\j,k\ne i}}
 |\delta J^{(i)}_{jk}|
 +\sum_{j\ne i}|\delta h^{(i)}_j|.
\]
Since spin and cosine factors have modulus at most one,
we add $\beta\sum_i\varepsilon_i$ to the log-partition
bound. Squared terms in $S(\gamma)$ are enclosed from above. For linear observable bounds, we similarly enclose
the difference $\delta h^{\mathrm{obs}}$ between the
stored tilted field and the exact shift
$h+(\mu/\beta)a$ in \eqref{eq:observable-shift},
adding $\beta\sum_j|\delta h^{\mathrm{obs}}_j|$
to the log-partition bound.

Rounding in $\beta=\gamma/(2c)$ may break the identity
$2c\beta=\gamma$ used in the elimination updates.
We therefore evaluate the envelope using the effective
curvature $c_{\mathrm{eff}}=\gamma/(2\beta)$.
Since the envelope offset decreases with curvature,
we use a lower enclosure of $c_{\mathrm{eff}}$ for
the offset and an upper enclosure of
$c_{\mathrm{eff}}\beta^2=\beta\gamma/2$ for the
quadratic coefficient.
Rounding in $\lambda=\beta\kappa/2$ changes
$\log F_\lambda$ by at most the error in $\lambda$,
since $|\partial_\lambda\log F_\lambda|\le1$.
We include this error in the upper bound.

\paragraph{Scalar envelopes.}
For the binary envelope, we enclose the positive solution
of $\tanh z=2cz$ by interval bisection and evaluate
$\log\cosh z-cz^2$ over the resulting interval to bound
$d(c)$ from above. For $c\ge1/2$, we use $d(c)=0$.

For the continuous envelope, to preserve an upper bound after truncation of the Fourier-Bessel function $F_\lambda$ \eqref{eq:fourier-bessel}, we explicitly
bound the omitted positive tail at order $K$. The power series
for $I_m$~\cite{DLMFBessel} gives
\begin{equation}
\begin{aligned}
 I_m(\lambda)
 =\sum_{j\ge0}\frac{q^{m+2j}}{j!(m+j)!}
 \le \frac{q^m}{m!}e^{q^2/(m+1)},\\
 \qquad q=\lambda/2.
 \end{aligned}
\end{equation}
Using $I_{2m}(r)\le I_0(r)$ and bounding successive
terms of $\sum_{m>K}q^m/m!$ by the ratio $q/(K+2)$,
the omitted positive tail satisfies
\begin{equation}
\begin{aligned}
     &T_K:=2\sum_{m>K}I_m(\lambda)I_{2m}(r)\\
 &\qquad \le
 \frac{2I_0(r)\,q^{K+1}e^{q^2/(K+2)}}
 {(K+1)!\,[1-q/(K+2)]},
 \\
 & \qquad\qquad\qquad\qquad 0\le q<K+2. 
\end{aligned}
\end{equation}
This tail upper bound is added to the truncated sum
before taking the logarithm. For $\lambda=0$,
$F_0(r)=I_0(r)$ is evaluated directly.
The maximum of $f(r)=\log F_\lambda(r)-cr^2$
lies in $[0,(2c)^{-1}]$. Adaptive subdivision uses
$f''\ge-2c$ and hence
\begin{equation}
 \sup_{r\in[a,b]}f(r)
 \le \max\{f(a),f(b)\}+\frac{c(b-a)^2}{4}
\end{equation}
to produce a global upper enclosure on the partition integral.

\paragraph{SDP and final bounds.}
Write the relaxation as minimising
$\operatorname{tr}(QY)$ over its feasible matrices $Y$. Each SDP minimisation problem has an associated
dual problem, whose feasible solutions provide lower
bounds on the original problem's optimal value.
For a numerical dual solution $y=(y_0,\ldots,y_N)$, define
\[
 D_{\rm B}(y)=\operatorname{diag}(y_0,y_1,\ldots,y_N)
\]
for binary spins, and
\[
 D_{\rm T}(y)=
 \operatorname{diag}(y_0,y_1,\ldots,y_N,y_1,\ldots,y_N)
\]
for continuous phases, with the ordering $(1,c,s)$.
A dual candidate provides a lower bound only if its slack
matrix $S=Q-D(y)$ is positive semidefinite.
Numerical solver output may violate this condition
slightly, so we bound the violation and correct the
reported objective.

For an adequate repair factor $\tau$ we compute an approximate Cholesky factor $R$ of
$S-\tau I$ and rigorously enclose the residual:
\[
 \|S-\tau I-R^\top R\|_\infty\le\rho.
\]
The enclosure accounts for all rounding errors in all
matrix operations involved. Here $\|\cdot\|_\infty$ is the maximum absolute row sum.
Because $R^\top R\succeq0$, this proves
$S\succeq(\tau-\rho)I$.
Consequently, for every feasible $Y\succeq0$,
\[
 \operatorname{tr}(SY)
 \ge(\tau-\rho)\operatorname{tr}Y.
\]
The normalisation constraints give
$\operatorname{tr}(D(y)Y)=\sum_{i=0}^N y_i$.
Using $Q=D(y)+S$ and $\operatorname{tr}Y=N+1$,
we therefore obtain
\[
 \operatorname{tr}(QY)
 \ge \sum_{i=0}^N y_i+(N+1)(\tau-\rho)
 =:L,
\]
as a certified lower energy bound, where we evaluate the expression using outward rounding toward a lower bound.

For magnetisation, we minimise $-O(x)=-a^\top x$
subject to $H(x)\le E$.
The energy constraint introduces a nonnegative dual
multiplier $\lambda$, giving the slack
\[
 S=Q_{-O}+\lambda(Q_H-Ee_0e_0^\top)-D_{\rm B}(y).
\]
Here $Q_H$ and $Q_{-O}$ represent $H(x)$ and $-O(x)$
under the lifting $Y=(1,x)(1,x)^\top$, and $e_0$ selects
the constant moment $Y_{00}=1$.
Applying the same correction to the dual slack as in the energy case gives
a certified lower bound on $-O(x)$. Its negative,
rounded upward, is the reported observable upper bound.

\section{Numerical settings and computational cost}
\label{app:numerical-settings}

We now elaborate on the numerical settings used in our computations and computational cost associated with our procedure.\paragraph{Parameter choices.}
We set $\gamma=g/s(J,h)$, where
\begin{align*}
 s(J,h)&=\max\!\left\{
 10^{-12},
 \sqrt{\frac{\|J\|_{\rm F}^2+\|h\|_2^2}{N}}
 \right\},
 \\
 g&\in\{0.08,0.15,0.30,0.60,0.90,1.20\}.
\end{align*}

For binary counts, we first evaluate the bound at ten equally
spaced values of the curvature $c$ between $0.06$ and $0.50$.
We also minimise the bound numerically over $c\in[0.04,0.50]$
using a bounded scalar optimiser, with an absolute tolerance
of $2\times10^{-7}$ in $c$ and a maximum of 80 iterations.
We retain the smallest bound obtained.

Continuous volumes use only
\[
 c\in\{0.06,0.08,0.11,0.15,0.20,0.25,0.30,0.40,0.50\},
\]
without local minimisation.
For each threshold, we retain the smallest enclosed
Chernoff exponent among the evaluated candidates,
including the trivial bound.

For magnetisation bounds, we evaluate the tilted models
$(J,h+ta)$ at
\[
 t\in\{0,0.15,0.30,0.60,1.0,1.5\},
\]
where $a$ is the magnetisation weight vector. For each tilt, we choose
\[
 \gamma=\frac{g}{s(J,h+ta)},
 \qquad g\in\{0.15,0.30,0.60,0.90\},
\]
using the scale function $s$ defined above. We evaluate eight equally spaced curvatures in $[0.08,0.50]$, without local refinement.

\medskip
\paragraph{SDP calculation.}
The first-order moment SDPs are constructed explicitly in
JuMP and solved using MOSEK, with one solver thread and
primal-feasibility, dual-feasibility, and relative-gap
tolerances of $10^{-8}$. No solver time limit is imposed.
The binary and full phase SDPs are solved separately.

\medskip
\paragraph{Computational cost.}
For $G$ profiles and dense coupling matrix, construction requires
$O(GN^3)$ time and $O(N^2)$ working memory.
Once constructed, the profiles are reused across different energy thresholds,
which require only scalar calculations.
In a representative Gaussian SK run at $N=1000$ (seed~1)
on an Apple MacBook Pro (M4 Pro, 24~GiB RAM),
using one thread each for BLAS and MOSEK,
profile construction took 240~s.
Additional binary and continuous threshold queries averaged
0.028~s and 0.00013~s, respectively.
The binary SDP calculation, including Julia startup,
took 10.9~s, followed by 0.42~s for residual validation.
The separate full phase SDP calculation took 30.9~s
using the already running Julia worker, followed by
4.81~s for residual validation.
Peak sampled process-tree memory across these stages
was approximately 2.75~GiB.
These are single representative measurements.

\bibliography{references}\end{document}